# Almost prime values of the order of abelian varieties over finite fields

Samuel Bloom[*]

March 9, 2018

**Abstract**

Let $E/\mathbb{Q}$ be an elliptic curve, and denote by $N(p)$ the number of $\mathbb{F}_p$-points of the reduction modulo $p$ of $E$. A conjecture of Koblitz, refined by Zywina, states that the number of primes $p \leq X$ at which $N(p)$ is also prime is asymptotic to $C_E \cdot X/\log(X)^2$, where $C_E$ is an arithmetically-defined non-negative constant. Following Miri-Murty (2001) and others, Y.R. Liu (2006) and David-Wu (2012) study the number of prime factors of $N(p)$. We generalize their arguments to abelian varieties $A/\mathbb{Q}$ whose adelic Galois representation has open image in $\operatorname{GSp}_{2g}\widehat{\mathbb{Z}}$. Our main result, after David-Wu, finds a conditional lower bound on the number of primes at which $\#A_p(\mathbb{F}_p)$ has few prime factors. We also present some experimental evidence in favor of a generalization of Koblitz's conjecture to this context.

## Contents

## 0 Introduction

Let $E/\mathbb{Q}$ be an elliptic curve without complex multiplication. Motivated by the heuristics of the Hardy-Littlewood Conjecture [HL23] and by cryptographic applications, a conjecture of Koblitz studies the number of points of $E_p$, the reduction of $E$ modulo a varying prime $p$. The conjecture posits:

**Conjecture 0.1** (Koblitz [Kob88], Conjecture A). *Suppose that every elliptic curve which is $\mathbb{Q}$-isogenous to $E$ (including $E$ itself) has trivial rational torsion. Then,*

$$\#\left\{p \leq x \text{ of good reduction} \,\Big|\, \#E_p\left(\mathbb{F}_p\right) \text{ is prime}\right\} \sim C_E \frac{x}{(\log x)^2}$$

*where $C_E$ is an explicit constant depending on the Galois representation of $E$.*

Koblitz's conjecture has been studied, refined, and generalized by many authors. In particular, the question of how often $\#E_p(\mathbb{F}_p)$ instead has *few prime factors* was first studied by Miri-Murty [MM01] with sieve-theoretic arguments. See Section 2 for more background and history. We study a generalization of this question to higher-dimensional abelian varieties, suggested in [Coj04; SW05] and others:

**Question 0.2.** *Let $A/\mathbb{Q}$ be an abelian variety of dimension $g$. At how many primes $p$ of good reduction does $\#A_p(\mathbb{F}_p)$ have few prime factors?*

---

[*] University of Maryland, College Park, MD, USA. *email:* bloom@math.umd.edu. *website:* www.math.umd.edu/~bloom

We briefly describe our main results, Theorems 1.3, 1.5, and 1.6. Let $\Omega(n)$ denote the number of prime factors of the integer $n$, counted with multiplicity, and $\omega(n)$ denote the number of prime factors of $n$ counted without multiplicity. Under the assumption that $A/\mathbb{Q}$ is in a certain sense "generic" and has no "congruence obstructions," and assuming a certain standard analytic hypothesis that generalizes the Riemann Hypothesis, we find

$$\#\left\{p \leq x \text{ of good reduction} \;\Big|\; \Omega\left(\#A_p\left(\mathbb{F}_p\right)\right) = r(g, \theta)\right\} \gg_A \frac{x}{(\log x)^2}$$

for an explicit function $r(g, \theta)$ which increases with $g$ but decreases with the strength of the analytic hypotheses, and

$$\#\left\{p \leq x \text{ of good reduction} \;\Big|\; \#A_p\left(\mathbb{F}_p\right) \text{ is prime}\right\} \ll_A \frac{x}{(\log x)^2},$$

with specified implicit constants. We show that $\omega\left(\#A_p(\mathbb{F}_p)\right)$ has *normal order* $\log \log p$ and, more precisely, that the quantity

$$\frac{\omega(\#A_p(\mathbb{F}_p)) - \log \log p}{\sqrt{\log \log p}}$$

follows a normal distribution as $p \to \infty$; this Erdös-Kac type result agrees with the heuristic that $\#A_p(\mathbb{F}_p)$ should act essentially like a random number which is approximately $p^g$. Our main innovation is the use of a proposition of [Cas+12] (labeled here as Proposition 4.10) and linear-algebraic arguments to compute the Chebotarev densities that appear when the sieve-theoretic arguments of [DW12] and [Liu06] are applied to higher-dimensional abelian varieties. We conclude by providing experimental evidence for the validity of Conjecture 6.1, an analogue of Conjecture 0.1 in the higher-dimensional setting, as suggested in [Coj04; SW05; Wen14].

## 0.1 Outline of paper.

We outline the structure of the paper. In Section 1, we provide precise statements of our results. In Section 2, we give background on Koblitz's conjecture, generalizations thereof, and the broader context of "questions of Lang-Trotter type" into which our results fit. In Section 3, we give preliminary information and known results regarding Galois representations, explicit Chebotarev Density theorems, and the weighted Greaves sieve that we will use in proving our results. In Section 4, we make preparations for the proof of our Theorems regarding prime-counting estimates and the counting of certain matrices. In Section 5, we prove Theorems 1.3, 1.5, and 1.6. In Section 6, we propose Conjecture 6.1 that generalizes Conjecture 0.1 and give experimental evidence for the Conjecture. In Section 7, we give concluding remarks.

## 0.2 Notations.

We use the standard Bachmann-Landau and Vinogradov notations for asymptotic growth of functions. In particular, subscripts on the asymptotic notations will denote dependency on the objects in the subscripts, so that an asymptotic notation without a subscript is absolute.

For a finite set $X$, we will write $\#X$ for the cardinality. We denote by $P_r$ the set of positive integers $n$ such that $\Omega(n) \leq r$. For a matrix $m$, denote by $\mathrm{char}_m(x)$ its characteristic polynomial.

We will use the letters $l$, $p$, $q$, and $\ell$ to denote rational prime numbers, $\mathfrak{p}$ to denote a prime ideal in a number field, and $\mathfrak{a}$ to denote an integral ideal in a number field. In a number field $L$, we will write $n_L$ or $n(L)$ for the degree of the extension $L/\mathbb{Q}$, and $d_L$ or $d(L)$ for the discriminant of the extension $L/\mathbb{Q}$.

For an abelian variety $A$ over a field $\kappa$, we will always use $\mathrm{End}(A)$ to denote the ring of endomorphisms of $A$ defined over the base field $\kappa$. For the sake of brevity, *we reserve $p$ and $\mathfrak{p}$ for places of $\kappa$ at which $A$ has good reduction.*

### Acknowledgments

I would like to thank my advisor, Prof. Larry Washington, for his advice and encouragement during this research. I would also like to thank Prof. Alina Cojocaru for kind and helpful communications.

# 1 Main Results

Throughout this article, we let $A/\mathbb{Q}$ be a principally polarized abelian variety of dimension $g \geq 1$ with conductor $N$ and adelic Galois representation $\widehat{\rho}\colon G_{\mathbb{Q}} \to \mathrm{GSp}_{2g}\,\widehat{\mathbb{Z}}$. We will call $A$ **generic**[1] if the image of $\widehat{\rho}$ is open in $\mathrm{GSp}_{2g}\,\widehat{\mathbb{Z}}$.

Recall that the Dedekind zeta function for a number field $L/\mathbb{Q}$,

$$\zeta_L(s) := \sum_{\mathfrak{a} \subset \mathcal{O}_L} \frac{1}{(\mathrm{N}\,\mathfrak{a})^s} = \prod_{\mathfrak{p} \subset \mathcal{O}_L} \left( \frac{1}{1 - (\mathrm{N}\,\mathfrak{p})^{-s}} \right)$$

has an analytic continuation to the entire complex plane, except for a simple pole at $s = 1$. Recall also that for a Galois extension $L/K$ of number fields, for each irreducible representation $\rho$ of $G := \mathrm{Gal}(L/K)$ we have the Artin $L$-function $L(s, \rho)$, that is in general known to be a meromorphic function on $\mathbb{C}$; moreover, we have the factorization

$$\zeta_L(s) = \zeta_K(s) \prod_{\substack{\rho \text{ non-triv. irred.} \\ \text{rep. of } G}} L(s, \rho)^{\deg(\rho)}$$

where $\deg(\rho)$ is the multiplicity of $\rho$ in the standard representation of $G$. Arithmetic information of $L$ and of $L/K$ is controlled by the zeros and coefficients of $\zeta_L$ and the Artin $L$-functions $L(s, \rho)$. In particular, there is the well-known

**Conjecture 1.1** (Artin's Holomorphy Conjecture (AHC) for $L/K$.)**.** *Let $\rho$ be a non-trivial irreducible representation of $\mathrm{Gal}(L/K)$. Then, $L(s, \rho)$ is holomorphic on $\mathbb{C}$.*

This conjecture is known for one-dimensional representations of $G$, since the Artin $L$-functions are then Hecke $L$-functions, which are known to be analytic on $\mathbb{C}$.

We will need to impose AHC as well as a generalization of the Riemann Hypothesis that asserts the existence of a zero-free half-plane region of $\zeta_L$ and of the $L$-functions for $L/K$. Ultimately, we will impose this hypothesis in Corollary 4.9, in the scenario that $L$ is a division field of $A$ and $K$ is a certain subfield.

**Hypothesis 1.2** ($\theta$-Hypothesis for $L/K$.)**.** *Let $1/2 \leq \theta < 1$, and $\mathcal{H}_\theta := \{s \in \mathbb{C} \mid \Re(s) > \theta\}$. Then, $\zeta_L(s)$ has no zeros in $\mathcal{H}_\theta$. Moreover, AHC holds for $L/K$, and the $L$-functions attached to irreducible representations of $\mathrm{Gal}(L/K)$ are zero-free on $\mathcal{H}_\theta$ as well.*

We require this analytic hypothesis to use the explicit error bound of [DW12], in the spirit of [MMS88], on the error terms in the Chebotarev Density Theorem. Following the argument of [DW12], we use these error bounds along with the weighted Greaves sieve (see Subsection 3.3) and find the following.

**Theorem 1.3.** *Suppose that $A$ is generic and that*

(Triv$_A$): *all of the abelian varieties over $\mathbb{Q}$ that are $\mathbb{Q}$-isogenous to $A$ have trivial rational torsion.*

*Assume the $\theta$-Hypothesis for the division fields of $A$ (i.e., for $\mathbb{Q}(A[n])/\mathbb{Q}$ for all $n$).[2] Then, for $x \gg_A 0$,*

$$\#\left\{ p \leq x \,\Big|\, \#A_p\left(\mathbb{F}_p\right) = P_r \right\} \geq B \cdot C_A \frac{x}{(\log x)^2}$$

*where $B$ is an explicit, absolute positive constant depending only on $g$, $C_A$ is an explicit non-negative constant depending on the Galois representation $\widehat{\rho}$ of $A$ (see (11)), and*

$$r = r(g, \theta) := \left\lceil \frac{(9/2)g^3 + (1/2)g}{1 - \theta} - \frac{1}{3} \right\rceil.$$

The utility of Theorem 1.3 is maximized once $\theta$ is small enough that $r(g, \theta) = r(g, 1/2) = 9g^3 + g$; thus, we obtain

---

[1] The terminology comes from the fact that for most dimensions $g$, the image of $\widehat{\rho}$ is open if the endomorphism ring of $A_{\overline{\mathbb{Q}}}$ is only $\mathbb{Z}$, which is true of "most" $A/\mathbb{Q}$. See, for instance, [Ser00b; Ser00a; Pin98] for more details.

[2] In fact, we only require the $\theta$-Hypothesis for $\mathbb{Q}(A[n])/\mathbb{Q}(A[n])^{B(n)}$, where $B(n)$ is a Borel subgroup of the Galois group of $\mathbb{Q}(A[n])/\mathbb{Q}$. We simplified the hypotheses here for the sake of readability.

**Corollary 1.4.** *Assume the hypotheses of Theorem 1.3, with*

$$\theta = 1 - \frac{(9/2)g^3 + (1/2)g}{9g^3 + g + 1/3}.$$

*Then, for $x \gg_A 0$,*

$$\#\left\{p \leq x \mid \#A_p(\mathbb{F}_p) = P_{9g^3+g}\right\} \geq B \cdot C_A \frac{x}{(\log x)^2}.$$

**Theorem 1.5.** *Suppose that $A$ is generic. Assume the $\theta$-Hypothesis for the division fields of $A$. Then, for all $\epsilon > 0$, for $x \gg_{A,\theta,\epsilon} 0$,*

$$\#\left\{p \leq x \mid \#A_p(\mathbb{F}_p) \text{ is prime}\right\} \leq \left(\frac{2g^2 + 3g + 6}{1 - \theta} + \epsilon\right) C_A \frac{x}{(\log x)^2}.$$

The constant $C_A$ is defined in (11) as an Euler product in terms of certain conjugacy classes attached to the Galois representation $\widehat{\rho}$. The assumption $(\text{Triv}_A)$ ensures that there is no "obvious" reason for all of the orders $\#A_p(\mathbb{F}_p)$ to share a common factor. We then understand there to be "congruence obstructions" to $\#A_p(\mathbb{F}_p)$ being prime infinitely often when $C_A = 0$. This possibility is the reason for the refinement by Zywina [Zyw11] of the constant $C_E$ for elliptic curves $E$.

We lastly follow the argument of Y.-R. Liu [Liu06], generalizing the Erdös-Kac Theorem, to show that $\#A_p(\mathbb{F}_p)$ essentially follows a normal distribution with *normal order* $\log \log p$.

**Theorem 1.6.** *Suppose that $A$ is generic. Assume the $\theta$-Hypothesis for the division fields of $A$ for some $\theta < 1$. Then, for all $\gamma \in \mathbb{R}$,*

$$\lim_{x \to \infty} \left(\frac{1}{\pi(x)} \#\left\{p \leq x \mid \frac{\omega(\#A_p(\mathbb{F}_p)) - \log \log p}{\sqrt{\log \log p}} \leq \gamma\right\}\right) = \frac{1}{\sqrt{2\pi}} \int_{-\infty}^{\gamma} e^{-t^2/2} \, dt.$$

## 2 Background.

Conjecture 0.1 fits within two broad families of questions of number-theoretic interest:

**Question 2.1.** *Given a "naturally-occurring" sequence $\mathcal{A}$ of integers (or tuples of integers), describe the subset $\Pi$ of terms which have a specified multiplicative-arithmetic behavior.*

**Question 2.2.** *Let $A$ be an abelian variety of dimension $g$ over a global field $L$. Let ♣ be a property of abelian varieties of dimension $g$ over finite fields. Describe*

$$\Pi = \Pi(A, \clubsuit) := \{\text{places } \mathfrak{p} : \mathcal{A}_\mathfrak{p} \text{ has } \clubsuit\},$$

*where $\mathcal{A}$ is the Néron model of $A$ over the appropriate one-dimensional scheme ($\text{Spec } \mathcal{O}_L$, resp. $C$, if $L$ is a number field, resp. the function field of a curve $C/\mathbb{F}_q$).*

By "describe" we mean either to give a "qualitative" description of $\Pi$ via congruence conditions, diophantine equations, and/or inequalities; or a "quantitative" description via an asymptotic estimate of the size of $\Pi$.

The family of Question 2.1 includes—among other questions—the Bateman-Horn Conjecture [BH62] which generalizes the Twin Prime Conjecture; the study of primes, pseudo-primes, and almost-primes in various intervals; and Artin's Conjecture on primitive roots modulo $p$ (see, for instance, [Mor12]). For an elliptic curve $E/\mathbb{Q}$, the family of Question 2.2 includes—among other questions—the Lang-Trotter Conjectures [LT76], the Sato-Tate Conjecture (see [Sut16] for an expository account), and the study of the structure of the group $E_p(\mathbb{F}_p)$ for a varying prime $p$ (see, for instance, [Coj04]). The generalizations of some these questions to higher-dimensional abelian varieties and/or over non-trivial number fields appear more difficult than their counterparts for elliptic curves over $\mathbb{Q}$. After Achter-Howe [AH17], we call this family "**questions of Lang-Trotter type.**"

The two families of Questions are very intertwined, and it would be hard to give a complete history. We content ourselves here with a brief history of Koblitz's conjecture and its generalizations. Koblitz based his conjecture on the heuristics behind the Hardy-Littlewood Conjecture: broadly, *unless there's an obstruction to it being otherwise, polynomials of degree d should (up to a correction factor that comes from congruence conditions) act like random number generators which, on input n, output a number on the order of $n^d$.*

From this heuristic, one finds a conjectural asymptotic count of the subset $\Pi$ by finding the expected value of a random variable for the probability distribution given by the heuristic. The heuristic probability distribution for Koblitz's Conjecture is based on the Sato-Tate distribution and the Galois representation of $E$, and states that

*unless there's an obstruction otherwise, the probability that $\#E_p(\mathbb{F}_p)$ is prime should be $C_E$ times the probability that a random number on the order of $p$ is prime.*

He thus conjectures that

$$\pi_E(x) := \#\left\{p \leq x \ \Big| \ \#E_p(\mathbb{F}_p) \text{ is prime}\right\} \sim \sum_{p \leq x} C_E \frac{1}{\log p}$$
$$\sim C_E \frac{x}{(\log x)^2}.$$

The motivation for Conjecture 0.1 was from cryptography: for the purposes of using the Elliptic Curve Discrete Logarithm Problem in a cryptographic protocol (for instance, in the Elliptic Curve Diffie-Helman key agreement protocol) one desires an elliptic curve over a large finite field having a prime number of points. Koblitz's suggestion was to choose an appropriate elliptic curve $E/\mathbb{Q}$, then reduce modulo appropriately large primes $p$ and find $\#E_p(\mathbb{F}_p)$ (via, for instance, the Schoof-Elkies-Atkin algorithm (see, e.g., [BSS99])) until $\#E_p(\mathbb{F}_p)$ is prime.

The question of finding a lower bound for $\pi_E(x)$ is still completely open: the author knows of no results showing even that $\pi_E(x) \to \infty$ for any specific elliptic curve. However, upper bounds are known. For $E/\mathbb{Q}$ *without* Complex Multiplication (**non-CM**), the first conditional and unconditional upper bounds are given by Cojocaru [Coj05] using the Selberg sieve. Zywina [Zyw08] improves upon these bounds by providing explicit asymptotic constants and extending the bounds to the case where $E$ is defined over a number field and may possibly have non-trivial torsion in its isogeny class. In the case of non-CM $E/\mathbb{Q}$, the best known conditional upper bound is given by David-Wu [DW12] who find

$$\pi_E(x) \leq \left(\frac{5}{1-\theta} + \epsilon\right) C_E \frac{x}{(\log x)^2}$$

for any $\epsilon > 0$, $x \gg_{\epsilon,\theta} 0$, assuming the $\theta$-Hypothesis for the division fields of $E$. For $E/\mathbb{Q}$ *with* CM, Cojocaru [Coj05] gives the unconditional upper bound $\pi_E(x) \ll_N x/(\log x)^2$.

Two approaches towards generalizing Conjecture 0.1 have yielded lower bounds. The first, which we do not pursue generalizing in this article, is to consider $\pi_E(x)$ *on average* for elliptic curves $Y^2 = X^3 + aX + b$ over $\mathbb{Q}$ in a family $\mathcal{C}(x)$, in the parameters $a$ and $b$ which vary in a rectangle that grows with $x$. That is, the approach is to consider the average

$$\lim_{x \to \infty} \left(\frac{1}{\#\mathcal{C}(x)} \sum_{E \in \mathcal{C}(x)} \pi_E(x)\right).$$

This was first considered in [BCD11] who show that the average is indeed $\sim \mathfrak{C}x/(\log x)^2$ if the rectangle for $\mathcal{C}(x)$ grows sufficiently quickly with respect to $x$, namely if $A, B > x^\epsilon$ and $AB > x(\log x)^{10}$. Here, $\mathfrak{C}$ is a positive constant to be thought of as an average of the $C_E$ for $E \in \mathcal{C}(x)$ as $x \to \infty$. They conclude then that "most" elliptic curves satisfy Conjecture 0.1; still, we cannot conclude Conjecture 0.1 for any specific curve. This result (and other "on average" results on the statistics of elliptic curves) has been improved; see, for instance, [DKS17].

The second approach, which we pursue in relation to abelian varieties, is to consider the question of *almost-prime* reductions of $E/\mathbb{Q}$. That is, this approach attempts to estimates

$$\pi_{E,r}(x) := \#\left\{p \leq x \ \Big| \ \#E_p(\mathbb{F}_p) \in P_r\right\}$$

for fixed $r$. This was first studied by Miri-Murty [MM01] who shows that for non-CM curves $E/\mathbb{Q}$ with trivial rational torsion, under GRH, $\pi_{E,16}(x) \gg x/(\log x)^2$. Steuding-Weng [SW05] improves this to $r = 9$ for non-CM curves, under GRH and the hypothesis (Triv$_E$). The best result for non-CM curves is by David-Wu [DW12], who show

$$\pi_{E,8}(x) \geq 2.778 \cdot C_E \frac{x}{(\log x)^2}$$

under the hypothesis (Triv$_E$) and the (11/21)-Hypothesis for the division fields of $E$. More precisely, their result is of the form

$$\pi_{E,r(\theta)}(x) \geq \frac{1.323}{1-\theta} C_E \frac{x}{(\log x)^2}$$

where the explicit function $r(\theta)$ decreases with the strength of the $\theta$-Hypothesis and is bounded below by 8. We will model our argument to theirs.

For elliptic curves over $\mathbb{Q}$ with CM, the situation is much better: Steuding-Weng first found $\pi_{E,3}(x) \gg x/(\log x)^2$ if $E$ is CM, under GRH and the hypothesis (Triv$_E$). Cojocaru [Coj05] improved this to $r = 5$ unconditionally. The best result for CM curves is by Iwaniec-Jiménez Urroz [IU10] and Jiménez Urroz [Jim08] who show unconditionally that

$$\#\{p \leq x \text{ of ordinary reduction} \mid \#E_p(\mathbb{F}_p) = d_E \cdot P_2\} \gg \frac{x}{(\log x)^2}$$

where $d_E = \gcd\{\#E_p(\mathbb{F}_p) \mid p \text{ of ordinary reduction}\}$.

More detailed statistical information of the function $p \mapsto \#E_p(\mathbb{F}_p)$ has been studied. In particular, Miri-Murty [MM01], Cojocaru [Coj05], and finally Y.-R. Liu [Liu06] find an Erdös-Kac result which provides a description of the "usual" behavior of $\#E_p(\mathbb{F}_p)$: they prove that, for any $\gamma \in \mathbb{R}$,

$$\lim_{x \to \infty} \left( \frac{1}{\pi(x)} \# \left\{ p \leq x \mid \frac{\omega(\#E_p(\mathbb{F}_p)) - \log \log p}{\sqrt{\log \log p}} \leq \gamma \right\} \right) = \frac{1}{\sqrt{2\pi}} \int_{-\infty}^{\gamma} e^{-t^2/2} \, \mathrm{d}t,$$

unconditionally if $E$ has CM, and conditionally on a $\theta$-Hypothesis on the division fields of $E$ if $E$ is non-CM. Liu concludes this normal distribution from a generalized version of the Erdös-Kac Theorem, which improves upon the generalized Hardy-Ramanujan result of Murty-Murty [MM84] that was used to study the coefficients of modular forms. We will use Liu's Theorem 3 to prove our Theorem 1.6.

Lastly, generalizations of Conjecture 0.1 to higher-dimensional abelian varieties have been suggested in and have begun to be studied. Weng [Wen14] computes the probability of the statement "$\ell \mid \#A(\mathbb{F}_p)$" for the reductions of a generic abelian variety of $\mathbb{Q}$, which we find in (17) in a different form. Weng [Wen15] and Spreckels [Spr17] also consider the "vertical" question of finding the probability, for fixed CM field $K$ and varying $p$, of the statement "$\exists A/\mathbb{F}_p$ with CM by $\mathcal{O}_K$ s.t. $\#A(\mathbb{F}_p)$ is prime," and conjecture an asymptotic behavior of a weighted counting function of such $p$, using the same heuristics as before.

## 3 Preliminaries.

### 3.1 Explicit Chebotarev Density Theorems

Let $L/\mathbb{Q}$ be a finite Galois extension with Galois group $G$ Let $C \subset G$ be a union of conjugacy classes. Denote by $\mathcal{P}(L/\mathbb{Q})$ the set of rational primes $p$ which ramify in $L/\mathbb{Q}$. Set

$$M(L/\mathbb{Q}) := n_L \prod_{\mathcal{P}(L/\mathbb{Q})} p.$$

Define the prime counting function for $C$,

$$\pi_C(x, L/\mathbb{Q}) := \# \{p \leq x : p \text{ unramified in } L/\mathbb{Q};\ \sigma_p \subseteq C\}$$

where $\sigma_p := \left(\frac{L/\mathbb{Q}}{p}\right)$ is the Artin symbol of $p$ in $L/\mathbb{Q}$. Recall that the **Chebotarev Density theorem** [Tsc26] states that as $x \to \infty$,

$$\pi_C(x, L/\mathbb{Q}) \sim \frac{\#C}{\#G} \int_2^x \frac{\mathrm{d}t}{\log t}.$$

We use the notation $\operatorname{li} x := \int_2^x \frac{dt}{\log t}$ for the logarithmic integral to $x$. We will use "explicit" versions of this theorem; that is, versions with bounds on the error term of the approximation. To state them, define the error term $R_C(x)$ via

$$\pi_C(x, L/\mathbb{Q}) = \frac{\#C}{\#G} \operatorname{li} x + R_C(x).$$

**Theorem 3.1** ([LO77; Ser81; MMS88; Mur97]). *Let the notation be as above.*

1. *Assume GRH for the Dedekind zeta function of $L/\mathbb{Q}$. Then,*

$$R_C(x) \ll (\#C) x^{1/2} \left( \frac{\log|d_L|}{n_L} + \log x \right)$$

2. *Assume GRH and AHC for $L/\mathbb{Q}$. Then,*

$$R_C(x) \ll (\#C)^{1/2} x^{1/2} \left( \log M(L/\mathbb{Q}) + \log x \right)$$

3. *Unconditionally, there exist positive constants $A, B, B'$ with $A$ effective and $B, B'$ absolute, such that if*

$$\log x \geq B'(\#G) \left( \log|d_L| \right)^2,$$

*then*

$$R_C(x) \ll \frac{\#C}{\#G} \operatorname{li} \left( x \cdot \exp\left( -B \frac{\log x}{\max\{|d_L|^{1/n_L}, \log|d_L|\}} \right) \right) + (\#\tilde{C}) x \cdot \exp\left( -A \sqrt{\frac{\log x}{n_L}} \right).$$

David-Wu extend the second statement to weaker assumptions. (In their notation, we set $K = \mathbb{Q}$.)

**Theorem 3.2** ([DW12]). *Let the notation be as above. Let $H \triangleleft G$ be a normal subgroup such that for all nontrivial irreducible characters $\rho$ of $\operatorname{Gal}(L^H/\mathbb{Q}) \cong G/H$, the Artin L-function $L(s, \rho)$ is holomorphic and is zero-free on the region $\{s \in \mathbb{C} \mid \Re(s) > \theta\}$. Suppose also that the product $HC \subseteq C$. Then,*

$$R_C(x) \ll \left( \frac{\#C}{\#H} \right)^{1/2} x^\theta n_L \left( \log M(L/\mathbb{Q}) + \log x \right).$$

This recovers the second part of Theorem 3.1 when $\theta = 1/2$ and $H$ is trivial.

We will also employ the following bound on the discriminant of $L/\mathbb{Q}$.

**Lemma 3.3** ([Ser81]). *Let the notation be as above. Then,*

$$\frac{n_L}{2} \sum_{\mathcal{P}(L/\mathbb{Q})} \log p \leq \log|d_L| \leq (n_L - 1) \sum_{\mathcal{P}(L/\mathbb{Q})} \log p + n_L \log n_L.$$

### 3.2 Galois Representations

Let $\kappa$ be a field and $\overline{\kappa}$ be an algebraic closure. Let $G_\kappa := \operatorname{Gal}(\overline{\kappa}/\kappa)$. Let $A/\kappa$ be a principally polarized abelian variety ("ppav") of dimension $g$ and conductor $N$. Recall that for any integer $n \geq 1$ not divisible by the characteristic of $\kappa$, the geometric torsion subgroup

$$A[n](\overline{\kappa}) \cong (\mathbb{Z}/n\mathbb{Z})^{2g}$$

is naturally a $G_\kappa$-module by action on the coordinates,

$$\rho_n : G_\kappa \to \operatorname{GSp}(A[n](\overline{\kappa}), e_n) \cong \operatorname{GSp}_{2g}(\mathbb{Z}/n\mathbb{Z}),$$

after choosing a symplectic basis of $A[n](\overline{\kappa})$ with respect to the Weil pairing $e_n$. We call $\rho_n$ the **mod-$n$ Galois representation** of $A$. We define the **$\ell$-adic Galois representation**, for $\ell \neq \operatorname{char} \kappa$, as the inverse limit

$$\rho_{\ell^\infty} := \varprojlim \rho_{\ell^n} : G_\kappa \to \operatorname{GSp}_{2g} \mathbb{Z}_\ell$$

after the symplectic identification

$$\varprojlim A[\ell^n](\overline{\kappa}) \cong \mathbb{Z}_\ell^{2g},$$

with the $\ell$-adic Weil pairing $e_{\ell^\infty}$ on the left and the standard symplectic pairing on the right. We then define the **adelic Galois representation**

$$\widehat{\rho} := \prod_\ell \rho_{\ell^\infty} : G_\kappa \to \prod_\ell \mathrm{GSp}_{2g} \mathbb{Z}_\ell \cong \mathrm{GSp}_{2g} \widehat{\mathbb{Z}} \tag{1}$$

The representations $\rho_n$, $\rho_{\ell^\infty}$, and $\widehat{\rho}$ are extremely important objects in the study of $A$. We will consider $\kappa = \mathbb{Q}$ and $\kappa$ a finite field.

**Notation 3.4.** *When we consider $\kappa = \mathbb{F}_p$, denote the Frobenius automorphism of $\overline{\mathbb{F}}_p/\mathbb{F}_p$ by $\mathrm{Frob}_p$. When we consider $\kappa = \mathbb{Q}$, for convenience we denote by $\mathrm{Frob}_p$ an absolute $p$-Frobenius automorphism, namely any choice of element in $\mathrm{Gal}(\overline{\mathbb{Q}}/\mathbb{Q})$ for which its image in $\mathrm{Gal}(L/\mathbb{Q})$, for any subextension $L$, has as its conjugacy class the Artin symbol $\left(\frac{L/\mathbb{Q}}{p}\right)$. It is well-known that $p$ is unramified in $\mathbb{Q}(A[l])/\mathbb{Q}$ when $\kappa = \mathbb{Q}$ (since $p \nmid lN$ under our notation), so that everything we will do is independent of this choice of conjugacy class.*

It is well-known that for $p \nmid N$ fixed and $\ell \neq p$ varying, the characteristic polynomial of Frobenius, $\mathrm{char}\,\rho_\ell(\mathrm{Frob}_p) \in \mathbb{Z}[x]$, is independent of $\ell$. We will thus without comment use the notation $\mathrm{char}(\mathrm{Frob}_p)$ or $\mathrm{char}_p$ for $\mathrm{char}(\rho_\ell(\mathrm{Frob}_p))$.

### 3.3 Simplified Greaves' Sieve

As in David-Wu, we use a simplified version of the weighted Greaves' Sieve for sieve problems of dimension 1, as given by Halberstam-Richert [HR85a; HR85b]. That is to say, in the notation of Halberstam-Richert, we will take $E = V$ and $T = U$.

For a set of primes $\mathcal{P}$, we use the notation

$$P(z) = \prod \left\{ p \,\Big|\, p \in \mathcal{P}, p \leq z \right\}.$$

For a list $\mathcal{A}$ of integers, and $d$ a positive integer, we use the notation

$$\mathcal{A}_d := \left\{ a \in \mathcal{A} \,\Big|\, a \equiv 0 \bmod d \right\}$$

**Theorem 3.5** (Simplified Greaves' Sieve, [HR85a; HR85b])**.** *Let $\mathcal{A}$ be a finite list of integers and $\mathcal{P}$ a set of primes such that the prime divisors of each $a \in \mathcal{A}$ are in $\mathcal{P}$. Let $y$ be a parameter, and $1/2 \leq U < 1$ and $V$ be constants such that*

$$V_0 \leq V \leq 1/4; \quad 1/2 \leq U < 1; \quad U + 3V \geq 1 \tag{2}$$

*where $V_0 = 0.074368\ldots$ is defined in [HR85a]. We suppose that there is a non-negative multiplicative function $w$ that satisfies the hypotheses*

$$w(p) = 0 \quad \text{for } p \notin \mathcal{P}, \tag{3}$$

$$0 < w(p) < p \quad \text{for } p \in \mathcal{P}, \tag{4}$$

$$\left| \sum_{\substack{p \in \mathcal{P}, \\ z_1 \leq p < z_2}} \frac{w(p)}{p} \log p - \log \frac{z_2}{z_1} \right| \leq A \quad \text{for } 2 \leq z_1 \leq z_2. \tag{5}$$

*Moreover, we suppose that there is an approximation $X \in \mathbb{R}^+$ to $\#\mathcal{A}$ and define the "remainders"*

$$r(\mathcal{A}, d) := \#\mathcal{A}_d - \frac{w(d)}{d} X \tag{6}$$

*for $d$ supported on $\mathcal{P}$. Define the sifting function*

$$H(\mathcal{A}, y^V, y^U) := \sum_{a \in \mathcal{A}} \gamma \left( \gcd(a, P(y^U)) \right)$$

*where*

$$\gamma(n) := \max\left\{0, \ 1 - \sum_{p|n, p \in \mathcal{P}} (1 - W(p))\right\},$$

*and where*

$$W(p) := \begin{cases} \dfrac{1}{U - V}\left(\dfrac{\log(p)}{\log(y)} - V\right) & \text{if } y^V \leq p < y^U, \\ 0 & \text{otherwise.} \end{cases}$$

*Then, we have the lower bound*

$$H(\mathcal{A}, y^V, y^U) \geq X \cdot V(y) \cdot \frac{2e^\gamma}{U - V}\left(J(U, V) + O\left(\frac{\log \log \log y}{(\log \log y)^{1/5}}\right)\right) \tag{7}$$

$$- (\log y)^{1/3} \left|\sum_{m < M, n < N, mn | P(y^U)} \alpha_m \beta_n \cdot r(\mathcal{A}, mn)\right|$$

*for any two real numbers $M, N$ such that*

$$MN = y; \quad M > y^U; \quad N > 1;$$

*with the $\alpha_m$ and $\beta_n$ certain real numbers in $[-1, 1]$; where*

$$V(y) := \prod_{p \leq y, p \in \mathcal{P}}\left(1 - \frac{w(p)}{p}\right);$$

$$J(U, V) := U \log \frac{1}{U} + (1 - U) \log \frac{1}{(1 - U)} - \log(4/3) + \alpha(V) - V \log 3 - V_0 \beta(V),$$

*where $\alpha(V)$ and $\beta(V)$ are certain non-negative numbers defined in [HR85b] as integrals, such that $\alpha(1/4) = \beta(1/4) = 0$.*

Halberstam-Richert apply this sieve to the problem of counting almost-primes in short intervals; see Theorem C from [HR85b]. Similarly, David-Wu apply the sieve to the problem of counting almost-prime orders of an elliptic curve $E/\mathbb{Q}$. They rely on the following Lemma (in a less general form), which uses the sifting function $H$ to detect these almost-prime orders. We will adapt this strategy to the higher-dimensional setting.

**Lemma 3.6** ([DW12]). *Let $\mathcal{A}$ be a finite list of positive integers, indexed by $\{p \leq x\}$, whose elements have all prime divisors in $\mathcal{P} = \{p \mid \gcd(p, M) = 1\}$. Suppose there exist real constants $U, V, \xi > 0$ and a positive integer $r$ such that $\max \mathcal{A} \leq (x^\xi)^{rU+V}$. (In the notation above, $y = x^\xi$.) Then,*

$$\#\left\{a \in \mathcal{A} \ \Big|\ \gcd(a, M) = 1; \ a = P_r\right\} \geq H\left(\mathcal{A}, (x^\xi)^V, (x^\xi)^U\right) - \sum_{(x^\xi)^V \leq p < (x^\xi)^U} \#\mathcal{A}_{p^2}.$$

### 3.4 Generalized Erdős-Kac Theorem

We will use the generalization of the Erdős-Kac Theorem by Y.R. Liu [Liu06] to prove Theorem 1.6. The classical Erdős-Kac Theorem states that the number of divisors of an integer $n$ has *normal order* $\log \log n$ and essentially follows a Gaussian distribution around that normal order.

**Theorem 3.7** (Erdös-Kac [EK40]).

$$\lim_{x \to \infty}\left(\frac{1}{x} \#\left\{n \leq x \ \Big|\ \frac{\omega(n) - \log \log n}{\sqrt{\log \log n}} \leq \gamma\right\}\right) = \frac{1}{\sqrt{2\pi}} \int_{-\infty}^{\gamma} e^{-t^2/2} \, dt.$$

Liu's generalization replaces $\omega(n)$ by $\omega(f(n))$ for functions $f$ with a particular shape. We state it here in the slightly more general form given by M. Xiong [Xio09]. In what follows, $S$ is an infinite subset of $\mathbb{N}$, and we use the notation $S(x) := \{n \in S \mid n \leq x\}$.

**Theorem 3.8** ([Liu06; Xio09]). *Suppose that $\#S(x^{1/2}) = o(\#S(x))$ as $x \to \infty$. Let $f : S \to \mathbb{N}$. For each prime $l$, choose functions $\lambda_l = \lambda_l(x)$ ("main term") and $e_l = e_l(x)$ ("error term") such that*

$$\frac{1}{\#S(x)} \# \left\{ n \in S(x) \,\middle|\, l \mid f(n) \right\} = \lambda_l + e_l.$$

*For increasing tuples $(l_1, \ldots, l_u)$ of distinct primes, we define functions $e_{l_1 \cdots l_u}(x)$ via*

$$\frac{1}{\#S(x)} \# \left\{ n \in S(x) \,\middle|\, l_1 \cdots l_u \mid f(n) \right\} = \left( \prod_{i=1}^{u} \lambda_{l_i} \right) + e_{l_1 \cdots l_u}.$$

*Suppose $\exists \beta \in (0,1], \exists c > 0$, independent of $x$, and a function $y = y(x)$ such that the following conditions hold:*

1. *for all $n \in S(x)$, the number of distinct prime divisors of $f(n)$ that are more than $x^\beta$ has a uniform upper bound (independent of $x$);*
2. $\sum_{y < l < x^\beta} \lambda_l = o(\sqrt{\log \log x})$;
3. $\sum_{y < l < x^\beta} |e_l| = o(\sqrt{\log \log x})$;
4. $\sum_{l < y} \lambda_l = c \log \log x + o(\sqrt{\log \log x})$;
5. $\sum_{l < y} \lambda_l^2 = o(\sqrt{\log \log x})$;
6. *for any $r \in \mathbb{N}$ and any integer $u$, $1 \leq u \leq r$,*

$$\sum_{\star} \left| e_{l_1 \cdots l_u}(x) \right| = o\left( (\log \log x)^{-r/2} \right)$$

   *where the sum $\sum_\star$ extends over all increasing tuples $(l_1, \ldots, l_u)$ of distinct primes $l_i < y(x)$.*

*Then, for $\gamma \in \mathbb{R}$,*

$$\lim_{x \to \infty} \left( \frac{1}{\#S(x)} \# \left\{ n \in S(x) \,\middle|\, \frac{\omega(f(n)) - c \log \log n}{\sqrt{\log \log n}} \leq \gamma \right\} \right) = \frac{1}{\sqrt{2\pi}} \int_{-\infty}^{\gamma} e^{-t^2/2} \, \mathrm{d}t.$$

## 4 Preparations for the Proof of Main Results.

Let $A/\mathbb{Q}$ be a generic abelian variety of conductor $N$. Recall that $p$ denotes a prime of good reduction for $A$, i.e., $p \nmid N$, and $l$ denotes a prime. Let

$$M = M_A := \prod \left\{ l \,\middle|\, \operatorname{im} \rho_{l^\infty} \neq \operatorname{GSp}_{2g} \mathbb{Z}_l \right\};$$

$$\mathcal{A} := \left\{ \#A_p(\mathbb{F}_p) \,\middle|\, p \leq x, \ \gcd(\#A_p(\mathbb{F}_p), M) = 1 \right\};$$

$$\mathcal{P} := \left\{ p \,\middle|\, p \nmid M \right\}.$$

Here, $\mathcal{A}$ is a *list,* i.e., might have repetition. We choose to omit from $\mathcal{A}$ those orders not coprime to $M$ so as to obtain the expected correction factor $C_A$ during the sieving process.

Our goal will be, assuming the $\theta$-Hypothesis, to show that for some choice of multiplicative function $w$, constants $U, V, \xi > 0$, and positive integer $r$, the hypotheses of Theorem 3.5 are satisfied; that

$$\text{right-hand side of (7)} \geq B \cdot C_A \frac{x}{(\log x)^2}, \tag{8}$$

for some constant $B > 0$; and that

$$\sum_{(x^\xi)^V \leq p < (x^\xi)^U} \#\mathcal{A}_{p^2} = o\left(\frac{x}{(\log x)^2}\right). \tag{9}$$

We will then choose such constants, depending on $\theta$, that minimize $r$. Theorem 1.3 will then follow from Lemma 3.6 with the constants we have chosen. After these computations, Theorem 1.5 will follow from the Selberg linear sieve, and Theorem 1.6 will follow from Theorem 3.8.

## 4.1 Divisibility of $\#A_p(\mathbb{F}_p)$

We recall some well-known facts about the Galois representations of $A$ and $A_p$. As in Subsection 3.2, for each $l$ we fix a $\mathbb{Z}_l$-basis of the $l$-adic Tate module of $A$ and of $A_p$ that is symplectic with respect to the Weil pairing. (For our purposes, we need not require any compatibility between these bases.) Thus, we may consider the $l$-adic Galois representations of $A$ and $A_p$ as taking values in $\operatorname{GSp}_{2g}(\mathbb{Z}_l)$.

Let $\pi_p \in \operatorname{End}(A_p)$ denote the Frobenius endomorphism. Recall the well-known theorem which states that for any abelian variety $B$ over a field $\kappa$, the restriction map $\operatorname{End}(B) \to \operatorname{End}_{\mathbb{Z}_l} T_l B$ is *injective*. Thus, we may consider $\pi_p$ as an element of $\operatorname{GSp}_{2g}(\mathbb{Z}_l)$.

**Theorem 4.1** (Weil Conjectures [Wei49; Gro66; Del73])**.** *The characteristic polynomial of $\pi_p \in \operatorname{GSp}_{2g} \mathbb{Z}_l$ has integer coefficients and is independent of $l$. Moreover, the eigenvalues of $\pi_p$ are $p$-**Weil numbers**. That is, all their embeddings into $\mathbb{C}$ have norm $\sqrt{p}$.*

Thus, the characteristic polynomial of $\pi_p$ has the form

$$\operatorname{char}_{\pi_p}(x) = x^{2g} + a_1 x^{2g-1} + \ldots + a_g x^g + p a_{g-1} x^{g-1} + p^2 a_{g-2} x^{g-2} + \ldots + p^g.$$

From now on, we consider Galois representations over $\kappa = \mathbb{Q}$. The following well-known lemma will allow us to detect information about $A_p$ from global information on $A$.

**Lemma 4.2.** *The conjugacy class of $\pi_p$ in $\operatorname{GSp}_{2g}(\mathbb{Z}_l)$ is $\rho_{l^\infty}(\operatorname{Frob}_p)$. In particular, $\operatorname{char}_{\pi_p} = \operatorname{char}_p$.*

From Lemma 4.2 and the observation that $\#A_p(\mathbb{F}_p) = \deg(\pi_p - \operatorname{id}_A) = \operatorname{char}_p(1)$, we immediately see that

**Lemma 4.3.** *For any $n \geq 1$, $n \mid \#A_p(\mathbb{F}_p) \iff \operatorname{char}_{\rho_n(\operatorname{Frob}_p)}(1) \equiv 0 \bmod n$.*

We are thus led to consider

**Definition 4.4.**

$$C(n) := \left\{ g \in \operatorname{Gal}(\mathbb{Q}(A[n])/\mathbb{Q}) \,\Big|\, \operatorname{char}_g(1) \equiv 0 \right\}, \tag{10}$$

so that, for $p \leq x$ such that $(\#A_p(\mathbb{F}_p), M) = 1$,

$$\#A_p(\mathbb{F}_p) \in \mathcal{A}_n \iff \rho_n(\operatorname{Frob}_p) \in C(n).$$

For convenience, for $n \geq 1$, set

**Notation 4.5.** $L_n := \mathbb{Q}(A[n])$, *and* $G(n) := \operatorname{Gal}(L_n/\mathbb{Q})$.

## 4.2 Setting up the sieve

We recall that the hypotheses of Theorem (3.5) require an approximation $X$ to $\#\mathcal{A}$ and a multiplicative function $w$ such that the "remainders" $r(\mathcal{A}, d)$ are small.

Mimicking the argument of David-Wu, we see that for squarefree $d$ that are supported on $\mathcal{P}$,

$$\#\mathcal{A}_d = \sum \left\{ 1 \;\middle|\; p \leq x, (\#A_p(\mathbb{F}_p), M_A) = 1, d \mid \#A_p(\mathbb{F}_p) \right\}$$

$$= \sum_{m \mid M_A} \mu(m) \cdot \sum \left\{ 1 \;\middle|\; p \leq x, dm \mid \#A_p(\mathbb{F}_p) \right\}$$

$$= \sum_{m \mid M_A} \mu(m) \cdot \pi_{C(dm)}(x, L_{dm}/\mathbb{Q})$$

$$= \sum_{m \mid M_A} \mu(m) \cdot \# \left\{ p \leq x \;\middle|\; \rho_d(\mathrm{Frob}_p) \subseteq C(d), \rho_m(\mathrm{Frob}_p) \subseteq C(m) \right\}$$

$$\sim \mathrm{li}(x) \frac{\#C(d)}{\#G(d)} \sum_{m \mid M_A} \left( \mu(m) \frac{\#C(m)}{\#G(m)} \right)$$

$$= \mathrm{li}(x) \frac{\#C(d)}{\#G(d)} \cdot \left( 1 - \frac{\#C'(M_A)}{\#G(M_A)} \right)$$

where

$$C'(M_A) := \left\{ g \in G(M_A) \;\middle|\; (\mathrm{char}_g(1), M_A) \neq 1 \right\}.$$

Thus, we choose

$$w(d) := \begin{cases} \frac{d \cdot \#C(d)}{\#G(d)} & d \text{ is supported on } \mathcal{P}; \\ 0 & \text{otherwise}; \end{cases} \qquad X := \mathrm{li}(x) \left( 1 - \frac{\#C'(M_A)}{\#G(M_A)} \right).$$

Then, $w$ is clearly multiplicative because of our assumption on $\widehat{\rho}$ and the Chinese Remainder Theorem. From these choices, the constant $C_A$ produced in the proof of 1.3 will then be

$$C_A = \left( 1 - \frac{\#C'(M_A)}{\#G(M_A)} \right) \lim_{y \to \infty} \left( \frac{V(y)}{\prod_{l < y}(1 - 1/l)} \right) = \frac{1 - \#C'(M_A)/\#G(M_A)}{\prod_{l \mid M_A}(1 - 1/l)} \prod_{l \nmid M_A} \frac{1 - \#C(l)/G(l)}{1 - 1/l}. \quad (11)$$

In order to show that $w$ satisfies the hypothesis (5), to find bounds on the remainders (6), and to show the bound (9), we will bound various Chebotarev densities, as well as find a bound on the size of $\#C(d)$. In the next subsection, before we begin computations, we describe a refinement of this argument, which we will use.

### 4.3 Exploiting subgroups of $\mathrm{GSp}_{2g}$

Using lemmas from [Ser81], David-Wu exploit the Borel and unipotent subgroups of $\mathrm{GL}_2$ and compare the prime counting functions for a Galois extension and a subextension to find the following.

**Theorem 4.6** ([DW12], Thm. 3.7). *Let $L/K$ be a Galois extension of number fields and $G = \mathrm{Gal}(L/K)$. Let $H \leq G$ and $C \subset G$ a union of conjugacy classes that intersects $H$. Let $C_H$ be the union of ($H$-)conjugacy classes in $H$ generated by $C \cap H$. Then,*

$$\pi_C(x, L/K) = \frac{|H|}{|G|} \frac{|C|}{|C_H|} \pi_{C_H}(x, L/L^H) + O\left( \frac{|C|}{|C_H| \cdot |G|} \log d_L + \frac{|H|}{|G|} \frac{|C|}{|C_H|} [L^H : \mathbb{Q}] x^{1/2} + [K : \mathbb{Q}] x^{1/2} \right).$$

Their idea (for $g = 1$) is to find subextensions

$$L_l^{H'(d)} \subset L_d^{H(d)} \subset L_d$$

where Theorem 4.6 applies to the extension $L_d^{H'(d)} \subset L_d$, and the second part of Theorem 3.1 applies to the subextension $L_d^{H'(d)} \subset L_d^{H(d)}$.

For this, they consider the subextensions

$$L_d^{B(d)} \subset L_d^{U(d)} \subset L_d$$

where $B(d)$ is the Borel subgroup of upper-triangular matrices in $\mathrm{GL}_2(\mathbb{Z}/d\mathbb{Z})$, and $U(d)$ is the subgroup of unipotent upper triangular matrices. Then,

$$B(d)/U(d) \cong \mathrm{Gal}\left(L_d^{U(d)}/L_d^{B(d)}\right)$$

is abelian, so AHC holds true in that extension. Thus, the second part of Theorem 3.1 applies to $L_d^{U(d)}/L_d^{B(d)}$.

We make preparations here to use the same idea in the setting of $g > 1$. *For the rest of this Section, we assume that $G(d) = \mathrm{GSp}_{2g}(\mathbb{Z}/d\mathbb{Z})$.*

**Notation 4.7.** *We set:*

- *$B(d)$ to be the (standard) Borel subgroup of $G(d)$, namely the subgroup of upper triangular matrices in $G(d)$;*
- *$U(d) \triangleleft B(d)$ to be the subgroup of unipotent matrices in $B(d)$;*
- *$C_B(d) := B(d) \cap C(d)$.*

We will also need to break up $G(d)$ into multiplicator cosets. Recall that for a commutative ring $R$ with unity,

$$\mathrm{GSp}_{2g}(R) := \left\{ M \in \mathrm{GL}_{2g}(R) \,\middle|\, \exists \mu \in R^\times \text{ s.t. } M^t J M = \mu J \right\}$$

where $J = \begin{pmatrix} 0 & I_g \\ -I_g & 0 \end{pmatrix}$ is the matrix for the standard symplectic form. We call the assignment $M \mapsto \mu$ the **multiplicator character** of $\mathrm{GSp}_{2g}$, and there is the exact sequence

$$1 \to \mathrm{Sp}_{2g}(R) \to \mathrm{GSp}_{2g}(R) \xrightarrow{\mu} R^\times \to 1. \tag{12}$$

For $m \in R^\times$, we define the $m$-**symplectic matrices,**

$$\mathrm{GSp}_{2g}^{(m)}(R) := \mu_R^{-1}(m)$$

and use the notation

$$G^{(m)}(d) := \mathrm{GSp}_{2g}^{(m)}(\mathbb{Z}/d\mathbb{Z}).$$

Now, we have the well-known

**Lemma 4.8.** *The characteristic polynomial of $M$ has the form*

$$\mathrm{char}_M(x) = x^{2g} + a_1 x^{2g-1} + \ldots + a_g x^g + m a_{g-1} x^{g-1} + m^2 a_{g-2} x^{g-2} + \ldots + m^g$$

*for some $a_i \in R$ and $m \in R^\times$.*

Thus, $B(d)/U(d)$ is the torus whose elements have coset representatives the diagonal matrices in $G(d)$ of the form $\begin{pmatrix} D & 0 \\ 0 & mD^{-1} \end{pmatrix}$ for a $g \times g$ invertible diagonal matrix $D$, so that $B(d)/U(d) \cong (\mathbb{G}_m(\mathbb{Z}/d\mathbb{Z}))^g \times \mathbb{G}_m(\mathbb{Z}/d\mathbb{Z})$. In particular, $B(d)/U(d)$ is abelian, so that, now in the context of $g \geq 1$, AHC holds true in the extension $L_d^{U(d)}/L_d^{B(d)}$. Therefore, by Theorem 3.2, we have

**Corollary 4.9.** *Assume the $\theta$-Hypothesis for the extensions $L_n/L_n^{B(n)}$. Then,*

$$\pi_{C_B(n)}(x, L_n/L_n^{B(n)}) = \frac{\#C_B(n)}{\#B(n)} \mathrm{li}(x) + R_n(x)$$

*where*

$$R_n(x) \ll \left(\frac{\#C_B(n)}{\#U(n)}\right)^{1/2} (\#B(n)) \cdot x^\theta \left(\log(M(L_n/L_n^{B(n)})) + \log x\right).$$

## 4.4 Fitting together the prime-counting estimates

In this subsection, we combine the discussion of Subsection 4.3 and the explicit Chebotarev Density Theorem. Using Theorem 4.6 with $G = G(n)$, $H = B(n)$, $C = C(n)$, $C_H = C_B(n)$, and $K = \mathbb{Q}$, we have

$$\pi_{C(n)}(x, L_n/\mathbb{Q}) = \frac{\#B(n)}{\#G(n)} \frac{\#C(n)}{\#C_B(n)} \pi_{C_B(n)}(x, L_n/L_n^{B(n)}) + Q_n(x)$$

where

$$Q_n(x) \ll \frac{\#C(n)}{\#C_B(n) \cdot \#G(n)} \log d_{L_n} + \frac{\#B(n)}{\#G(n)} \frac{\#C(n)}{\#C_B(n)} [L_n^{B(n)} : \mathbb{Q}] x^{1/2} + x^{1/2}. \tag{13}$$

Plugging Corollary 4.9 into the above and canceling factors, we have

$$\pi_{C(n)}(x, L_n/\mathbb{Q}) = \frac{\#C(n)}{\#G(n)} \mathrm{li}(x) + \frac{\#B(n)}{\#G(n)} \frac{\#C(n)}{\#C_B(n)} R_n(x) + Q_n(x) \tag{14}$$

where $R_n$ and $Q_n$ have their respective bounds as above (with the bound on $R_n$ assuming the $\theta$-Hypothesis).

## 4.5 Counting matrices

In this subsection, we compute and gather estimates on the sizes of the subsets of $G(d)$ (and their ratios) that have appeared. *We maintain the assumption that $G(d) = \mathrm{GSp}_{2g}(\mathbb{Z}/d\mathbb{Z})$.*

To begin, there is the well known formula

$$\# \mathrm{Sp}_{2g} \mathbb{F}_l = l^{g^2} \prod_{i=1}^{g} \left( l^{2i} - 1 \right) = l^{2g^2+g} - l^{2g^2+g-2} + O_g\left( l^{2g^2+g-6} \right). \tag{15}$$

From this and the exact sequence (12), we have

$$\#G(l) = (l-1) l^{g^2} \prod_{i=1}^{g} \left( l^{2i} - 1 \right) = l^{2g^2+g+1} - l^{2g^2+g} + O_g\left( l^{2g^2+g-1} \right). \tag{16}$$

Recall Definition (10). For convenience, for an integer $m$, denote

$$C^{(m)}(d) := C(d) \cap G^{(m)}(d).$$

From Castryck et al. [Cas+12], we have

**Proposition 4.10.**

$$\frac{\#C^{(m)}(l)}{\#G^{(m)}(l)} = \begin{cases} -\sum_{r=1}^{g} l^r \prod_{j=1}^{r} (1-l^{2j})^{-1} & \text{if } l \mid m-1, \\ -\sum_{r=1}^{g} \prod_{j=1}^{r} (1-l^j)^{-1} & \text{otherwise.} \end{cases}$$

Thus,

$$\#C(l) = \sum_{m \in (\mathbb{Z}/l\mathbb{Z})^\times} \frac{\#C^{(m)}(l)}{\#G^{(m)}(l)} \#G^{(m)}(l)$$

$$= \sum_{m \in (\mathbb{Z}/l\mathbb{Z})^\times} \frac{\#C^{(m)}(l)}{\#G^{(m)}(l)} \cdot \frac{\#G(l)}{l-1}$$

$$= \frac{\#G(l)}{l-1} \cdot \left( -\sum_{r=1}^{g} l^r \prod_{j=1}^{r} (1-l^{2j})^{-1} + (l-2)\left( -\sum_{r=1}^{g} \prod_{j=1}^{r} (1-l^j)^{-1} \right) \right) \tag{17}$$

$$= l^{2g^2+g} - 3 l^{2g^2+g-2} + O_g(l^{2g^2+g-3}) \tag{18}$$

When $g = 1$, (17) yields $\frac{\#C(l)}{\#G(l)} = \frac{l^2-2}{(l-1)(l^2-1)}$, which agrees with the density written in David-Wu.

Next, we count $\#B(l)$. Since

$$B(l) = \left\{ M \in \mathrm{GL}_{2g}\,\mathbb{Z}/l\mathbb{Z} \,\Big|\, M \text{ upper triangular}, M \in \mathrm{GSp}_{2g}\,\mathbb{Z}/l\mathbb{Z} \right\},$$

then $B(l)$ consists of $M = \begin{pmatrix} T_1 & A \\ 0 & T_2 \end{pmatrix}$ with $T_i$ upper-triangular, such that for some $\mu$,

$$\begin{pmatrix} T_1 & A \\ 0 & T_2 \end{pmatrix}^t J \begin{pmatrix} T_1 & A \\ 0 & T_2 \end{pmatrix} = \mu J,$$

i.e.

$$T_1^t T_2 = \mu I, \qquad A^t T_2 = T_2^t A;$$

i.e.

$$T_2 = \mu \left(T_1^t\right)^{-1}; \qquad A = \mu^{-1} T_1 R$$

for some symmetric matrix $R$. That is,

$$B(l) = \left\{ \begin{pmatrix} T & \mu^{-1} R T^t \\ 0 & \mu^{-1}(T^t)^{-1} \end{pmatrix} \right\} \tag{19}$$

and thus

$$\#B(l) = (l-1)\left((l-1)^g \cdot l^{g(g-1)/2}\right)\left(l^{g(g+1)/2}\right)$$
$$= (l-1)^{g+1} \cdot l^{g^2},$$

and $\#U(l) = l^{g^2}$. From this description we also see that

$$\frac{\#C_B(l)}{\#B(l)} = 1 - \frac{\#\left\{(T,\mu) \,\Big|\, T \text{ does not have 1 as an eig.val.}, \mu^{-1} \notin \{\text{eig.vals. of } T\}\right\}}{\#\{(T,\mu)\}}$$

so that

$$\frac{\#C_B(l)}{\#B(l)} \leq 1 - \frac{\#\left\{T \,\Big|\, T \text{ does not have 1 as an eig.val.}\right\} \cdot (l-1-g)}{(l-1)^{g+1} \cdot l^{g(g-1)/2}};$$

$$\frac{\#C_B(l)}{\#B(l)} \geq 1 - \frac{\#\left\{T \,\Big|\, T \text{ does not have 1 as an eig.val.}\right\} \cdot (l-2)}{(l-1)^{g+1} \cdot l^{g(g-1)/2}}$$

and thus

$$1 - \frac{(l-2)^g(l-2)}{(l-1)^{g+1}} \leq \frac{\#C_B(l)}{\#B(l)} \leq 1 - \frac{(l-2)^g(l-1-g)}{(l-1)^{g+1}} \tag{20}$$

so that $\#C_B(l)/\#B(l) \asymp_g 1/l$.

We also record here for future use that, by the same reasoning,

$$1 - \frac{(l^2-l-1)^g(l^2-l-2)}{(l^2-l)^{g+1}} \leq \frac{\#C_B(l^2)}{\#B(l^2)} \leq 1 - \frac{(l^2-l-1)^g(l^2-l-g)}{(l^2-l)^{g+1}}$$

so that $\#C_B(l^2)/\#B(l^2) \asymp_g 1/l$.

Next, $\#C(l^2)$. It will suffice for our purposes to have an upper bound on $\frac{\#C(l^2)}{\#G(l^2)}$.

**Lemma 4.11.** $\frac{\#C(l^2)}{\#G(l^2)} = O_g\left(\frac{1}{l^2}\right).$

*Proof.* We write

$$\frac{\#C(l^2)}{\#G(l^2)} = \frac{\#C(l^2)}{\#C(l)} \cdot \frac{\#C(l)}{\#G(l)} \cdot \frac{\#G(l)}{\#G(l^2)}$$

Now, consider the mod-$l$ reduction map, which is surjective by Hensel's Lemma (see, e.g., pg. 177 of [Mum99]):

$$1 \to K \to \operatorname{GSp}_{2g} \mathbb{Z}/l^2\mathbb{Z} \xrightarrow{\phi_l} \operatorname{GSp}_{2g} \mathbb{Z}/l\mathbb{Z} \to 1,$$

where

$$K = \left(I + l \cdot M_{g \times g}(\mathbb{Z}/l^2\mathbb{Z})\right) \cap \operatorname{GSp}_{2g} \mathbb{Z}/l^2\mathbb{Z}.$$

Then, $\#G(l)/\#G(l^2) = 1/\#K$. From earlier discussion, we also have $\#C(l)/\#G(l) = O_g(1/l)$.

It remains to bound $\frac{\#C(l^2)}{\#C(l)}$. Note that $C(l^2) \subset \phi_l^{-1}(C(l))$, so in particular the product $K \cdot C(l^2) \subset \phi_l^{-1}(C(l))$. We show that

$$\#C(l^2) \leq \frac{2g}{l} \# \left(K \cdot C(l^2)\right) \leq \frac{2g}{l} \#K \cdot \#C(l);$$

the second inequality is obvious.

Consider the subgroup of scalar matrices $S := \left((1 + l\mathbb{Z})/l^2\mathbb{Z}\right) \cdot I \subset K$. For $\alpha I \in S$ and $M \in C(l^2)$, the product $\alpha M$ is in $C(l^2)$ only when one of the eigenvalues, say $\beta$, of $M$ is such that $\alpha\beta \equiv 1 \mod l^2$. But since $\alpha \in (1 + l\mathbb{Z})/l^2\mathbb{Z}$, the equation $\alpha\beta \equiv 1 \mod l^2$ has only one solution $\beta$.

Thus, accounting for the possible multiplicity of the eigenvalues of $M$, we have

$$\#S\{M\} \cap C(l^2) \leq 2g.$$

So, partition $C(l^2)$ into subsets of orbits under $S$; that is, form the set

$$C(l^2)/S := \left\{ SM \cap C(l^2) \mid M \in C(l^2) \right\}.$$

Then, $\#C(l^2)/S \geq \#C(l^2)/2g$, so that

$$\# \left(S \cdot C(l^2)\right) = \#S \cdot \# \left(C(l^2)/S\right) \geq \frac{l}{2g} \#C(l^2).$$

But certainly $\# \left(K \cdot C(l^2)\right) \geq \# \left(S \cdot C(l^2)\right)$, and thus the desired inequality follows. $\square$

Lastly, we record the following formulas. A short proof of the first formula is given in `https://mathoverflow.net/questions/87904`. We will only use this formula in the case $k = 2$. The proof of the second formula is clear from (19).

**Lemma 4.12.**

$$\#G(l^k) = (l-1)l^{(2k-1)g^2 + (k-1)g + 1} \prod_{i=1}^{g}(l^{2i} - 1). \qquad \#B(l^2) = (l-1)^{g+1} l^{2g^2 + g + 1}.$$

### 4.6 Verifying the sieve hypothesis (5)

We now verify hypothesis (5). Recall that we defined

$$w(d) := \begin{cases} \frac{d \cdot \#C(d)}{\#G(d)} & d \text{ is supported on } \mathcal{P}; \\ 0 & \text{otherwise}; \end{cases}$$

.

From (17), we have
$$\frac{w(l)}{l} = \frac{1}{l} + O_g\left(\frac{1}{l^2}\right)$$
so that for $z_1 < z_2$,
$$\left|\sum_{z_1 \leq l < z_2,\ l \in \mathcal{P}} \frac{w(l)}{l} \log l - \log \frac{z_2}{z_1}\right| = \left|\sum_{z_1 \leq l < z_2,\ l \in \mathcal{P}} \left(\frac{1}{l} + O_g\left(\frac{1}{l^2}\right)\right) \log l - \log \frac{z_2}{z_1}\right|$$
$$= \left|\sum_{z_1 \leq l < z_2,\ l \in \mathcal{P}} \frac{1}{l} \log l - \log \frac{z_2}{z_1}\right| + \left|\sum_{l \in \mathcal{P}} O_g\left(\frac{1}{l^2}\right) \log l\right|.$$

By the comparison test for series, the second term is $O_g(1)$. Recall now one of Mertens' theorems,

**Theorem 4.13** ([Mer74]). *For all $n \geq 2$,*
$$\left|\sum_{l \leq n} \frac{\log(l)}{l} - \log(n)\right| \leq 2$$

Thus, via the Triangle Inequality, hypothesis (5) is verified, so that Theorem (3.5) applies, and so for any valid choice of constants, the lower bound (7) holds.

## 5 Proof of Main Results.

We now combine the estimates of this section and the theorems of Subsection 3.1 in order to show the existence of constants $U, V, \xi, r$ that guarantee the lower bound (8) and the upper bound (9). We will box the constraints on the constants as we determine them.

First, the hypothesis of Lemma 3.6 that $\max \mathcal{A} \leq (x^\xi)^{rU+V}$ requires, by earlier discussion, that $\boxed{g < \xi(rU + V).}$

### 5.1 Ensuring (8)

We begin with the lower bound (8). Recall that we wish to show that

$$X \cdot V(y) \cdot \frac{2e^\gamma}{U - V} \left(J(U, V) + O\left(\frac{\log \log \log y}{(\log \log y)^{1/5}}\right)\right)$$
$$- (\log y)^{1/3} \left|\sum_{m < M, n < N, mn | P(y^U)} \alpha_m \beta_n \cdot r(\mathcal{A}, mn)\right| \geq B \cdot C_A \cdot \frac{x}{(\log x)^2}$$

where
$$X := \mathrm{li}(x)\left(1 - \frac{\#C'(M_A)}{\#G(M_A)}\right)$$
$$V(y) := \prod_{p \leq y, p \in \mathcal{P}}\left(1 - \frac{w(p)}{p}\right),$$
$$C_A := \frac{1 - \#C'(M_A)/\#G(M_A)}{\prod_{l | M_A}(1 - 1/l)} \prod_{l \nmid M_A}\left(\frac{1 - \#C(l)/\#G(l)}{1 - 1/l}\right),$$
$$r(\mathcal{A}, d) := \#\mathcal{A}_d - \frac{w(d)}{d} \cdot \left(1 - \frac{\#C'(M_A)}{\#G(M_A)}\right) \mathrm{li}(x),$$
$$w(d) := \begin{cases} \frac{d \cdot \#C(d)}{\#G(d)} & d \text{ is supported on } \mathcal{P}; \\ 0 & \text{otherwise;} \end{cases}$$

and $M, N, \alpha_m, \beta_n, \alpha(V), \beta(V)$ are as in previous notation.

As in Lemma 3.6, we choose
$$y = x^\xi.$$

Now, following the argument of David-Wu, assuming that $x^\xi > M_A$,

$$V(x^\xi) = \prod_{l<x^\xi, l \nmid M_A} \left(1 - \frac{1}{l}\right) \prod_{l<x^\xi, l \nmid M_A} \left(\frac{1 - \#C(l)/\#G(l)}{1 - 1/l}\right)$$

$$= \prod_{l<x^\xi} \left(1 - \frac{1}{l}\right) \prod_{l|M_A} \left(1 - \frac{1}{l}\right)^{-1} \prod_{l<x^\xi, l \nmid M_A} \left(\frac{1 - \#C(l)/\#G(l)}{1 - 1/l}\right)$$

$$\overset{\text{Mertens}}{\sim} \frac{e^{-\gamma}}{\xi \log x} \cdot C_A \cdot \left(1 - \#C'(M_A)/\#G(M_A)\right)^{-1} \cdot \prod_{l>x^\xi} \left(\frac{1 - \#C(l)/\#G(l)}{1 - 1/l}\right)^{-1}$$

where the asymptotic $\sim$ is as $x^\xi \to \infty$.

Then, considering the "remainder" $r(\mathcal{A}, d)$ for squarefree $d$ supported on $\mathcal{P}$, we have

$$r(\mathcal{A}, d) = \sum_{m|M_A} \left(\mu(m) \cdot \pi_{C(dm)}(x, L_{dm}/\mathbb{Q})\right) - \frac{\#C(d)}{\#G(d)} \left(1 - \frac{\#C'(M_A)}{\#G(M_A)}\right) \mathrm{li}(x).$$

But since $G(dm) = G(d) \times G(m)$ for $m \mid M_A$ and $d$ supported outside of $M_A$, then, using (14)

$$\sum_{m|M_A} \mu(m) \cdot \pi_{C(dm)}(x, L_{dm}/\mathbb{Q}) =$$

$$= \sum_{m|M_A} \mu(m) \cdot \left(\frac{\#C(dm)}{\#G(dm)} \mathrm{li}(x) + \frac{\#B(dm)}{\#G(dm)} \frac{\#C(dm)}{\#C_B(dm)} R_{dm}(x) + Q_{dm}(x)\right)$$

$$= \sum_{m|M_A} \mu(m) \cdot \left(\frac{\#C(d)}{\#G(d)} \frac{\#C(m)}{\#G(m)} \mathrm{li}(x) + \frac{\#B(dm)}{\#G(dm)} \frac{\#C(dm)}{\#C_B(dm)} R_{dm}(x) + Q_{dm}(x)\right)$$

$$= \left(1 - \frac{\#C'(M)}{\#G(M)}\right) \frac{\#C(d)}{\#G(d)} \mathrm{li}(x)$$

$$+ \frac{\#B(d)}{\#G(d)} \frac{\#C(d)}{\#C_B(d)} \sum_{m|M_A} \mu(m) \frac{\#B(m)}{\#G(m)} \frac{\#C(m)}{\#C_B(m)} R_{dm}(x) + \sum_{m|M_A} \mu(m) Q_{dm}(x).$$

But $R_{dm}(x) \ll_A R_d(x)$ and $Q_{dm}(x) \ll_A Q_d(x)$, so we have

$$r(\mathcal{A}, d) = \frac{\#B(d)}{\#G(d)} \frac{\#C(d)}{\#C_B(d)} \cdot O_A(R_d(x)) + O_A(Q_d(x)).$$

By the Chinese Remainder Theorem, (17), and (20),

$$\frac{\#B(d)}{\#G(d)} \frac{\#C(d)}{\#C_B(d)} \leq \prod_{l|d} \left(1 - \frac{(l-2)^g(l-1-g)}{(l-1)^{g+1}}\right)^{-1} \left(\frac{1}{l-1} \cdot \left(-\sum_{r=1}^{g} l^r \prod_{j=1}^{r}(1-l^{2j})^{-1} + (l-2)\left(-\sum_{r=1}^{g} \prod_{j=1}^{r}(1-l^j)^{-1}\right)\right)\right)$$

$$\ll_g \prod_{l|d}(l-1) \cdot \frac{1}{l-1} \cdot \frac{l^2-1}{(l^2-2)}$$

$$< \prod_l \left(1 + \frac{1}{l^2-2}\right) = O(1).$$

Next, from Corollary 4.9 we have

$$R_d(x) \ll \left(\frac{\#C_B(d)}{\#U(d)}\right)^{1/2} (\#B(d)) \cdot x^\theta \left(\log(M(L_d/L_d^{B(d)})) + \log x\right)$$

and so, again by (20) and the Chinese Remainder Theorem,

$$R_d(x) \ll \left(\left(\frac{1}{\#U(d)}\right)^{1/2}\left(\frac{C_B(d)}{\#B(d)}\right)^{1/2}(\#B(d))^{3/2}\right)x^\theta\left(\log(M(L_d/L_d^{B(d)}))+\log x\right)$$

$$\ll_g d^{g^2+(3/2)g+1}x^\theta(\log(d)+\log(x)).$$

Next, from (13), using Lemma 3.3,

$$Q_d(x) \ll \frac{\#C(d)}{\#C_B(d)\cdot\#G(d)}\log d_{L_d} + \frac{\#B(d)}{\#G(d)}\frac{\#C(d)}{\#C_B(d)}[L_d^{B(d)}:\mathbb{Q}]x^{1/2}+x^{1/2}$$

$$\ll \frac{\#C(d)}{\#C_B(d)\cdot\#G(d)}\log d_{L_d} + \frac{\#C(d)}{\#C_B(d)}x^{1/2}+x^{1/2}$$

$$\ll \frac{1-\epsilon}{2g\cdot\#B(d)}\log d_{L_d} + \frac{\#G(d)}{\#B(d)}\frac{1}{2g}(1-\epsilon)x^{1/2}+x^{1/2}$$

$$\ll_A \frac{1}{d^{g+1}d^{g^2}}\cdot d^{2g^2+g+1}\log(d) + \frac{d^{2g^2+g+1}}{d^{g+1}d^{g^2}}x^{1/2}$$

$$= d^{g^2}\left(x^{1/2}+\log(d)\right)$$

Thus, since $\theta \geq 1/2$,

$$r(\mathcal{A},d) \ll_g d^{g^2+(3/2)g+1}x^\theta(\log(d)+\log(x)) + d^{g^2}\left(x^{1/2}+\log(d)\right)$$

$$\ll_\epsilon d^{g^2+(3/2)g+1}x^{\theta+\epsilon}$$

Now, by the Triangle Inequality, the sum

$$\left|\sum_{m<M,n<N,mn|P(y^U)}\alpha_m\beta_n\cdot r(\mathcal{A},mn)\right| \leq \sum_{m<M,n<N,mn|P(y^U)}|r(\mathcal{A},mn)|$$

$$\ll_{A,\epsilon} x^{\theta+\epsilon}\log(x)\sum_{m<M,n<N,mn|P(x^{\xi U})}(mn)^{g^2+(3/2)g+1}$$

Since $P(x^{\xi U})$ is squarefree, we note that for any non-negative function $f(t)$, since $U<1$, we have

$$\sum_{m<M,n<N,mn|P(x^{\xi U})}f(mn) \leq \sum_{d\leq x^\xi}\mu(d)^2 3^{\omega(d)}f(d).$$

But $3^{\omega(d)}\leq(3/2)d$, and, of course, $\mu(d)^2\leq 1$. Thus, integrating by parts, the sum above is

$$\ll x^{\theta+\epsilon}\int_1^{x^\xi} d^{g^2+(3/2)g+2}\cdot\mathrm{d}(\text{sq.free. ints.})$$

$$\ll x^{\theta+\epsilon+\xi(g^2+(3/2)g+3)}.$$

Thus, finally, the lower bound (8) will be satisfied if

$$\boxed{\theta+\epsilon+\xi\left(g^2+\frac{3}{2}g+3\right)<1}$$

with the constant

$$B = \xi^{-1}\cdot\frac{J(U,V)}{U-V}.$$

## 5.2 Ensuring (9)

We now ensure the lower bound (9), namely that

$$\sum_{(x^\xi)^V \leq l < (x^\xi)^U} \#\mathcal{A}_{l^2} = o\left(\frac{x}{(\log x)^2}\right).$$

We have that

$$\#\mathcal{A}_{l^2} = \frac{\#C(l^2)}{\#G(l^2)} \operatorname{li}(x) + \frac{\#B(l^2)}{\#G(l^2)} \frac{\#C(l^2)}{\#C_B(l^2)} R_{l^2}(x) + Q_{l^2}(x)$$

$$= O_g\left(\frac{1}{l^2}\right) \operatorname{li}(x) + O_g\left(\frac{1}{l^2}\right) \cdot O_g(l) R_{l^2}(x) + Q_{l^2}(x)$$

$$= O_g\left(\frac{1}{l^2}\right) \operatorname{li}(x) + O_g\left(\frac{1}{l}\right) R_{l^2}(x) + Q_{l^2}(x)$$

where

$$R_{l^2}(x) \ll \left(\frac{\#C_B(l^2)}{\#U(l^2)}\right)^{1/2} (\#B(l^2)) \cdot x^\theta \left(\log(M(L_{l^2}/L_{l^2}^{B(l^2)})) + \log x\right)$$

$$\ll_g \left(\left(\frac{1}{\#U(l^2)}\right)^{1/2} \left(\frac{\#C_B(l^2)}{\#B(l^2)}\right)^{1/2} (\#B(l^2))^{3/2}\right) \cdot x^\theta (\log l + \log x)$$

$$\ll_g \left(\left(\frac{1}{l^{3g^2+1}}\right)^{1/2} \left(\frac{1}{l}\right)^{1/2} ((l-1)^{g+1} l^{4g^2-g})^{3/2}\right) \cdot x^\theta \log x$$

$$\ll l^{(9/2)g^2 + 1/2} x^\theta \log x$$

and

$$Q_{l^2}(x) \ll \frac{\#C(l^2)}{\#C_B(l^2) \cdot \#G(l^2)} \log d_{L_{l^2}} + \frac{\#B(l^2)}{\#G(l^2)} \frac{\#C(l^2)}{\#C_B(l^2)} [L_{l^2}^{B(l^2)} : \mathbb{Q}] x^{1/2} + x^{1/2}$$

$$\ll \frac{1}{\#B(l^2)} \cdot O_g\left(\frac{1}{l}\right) \left(\#G(l^2) \log(\#G(l^2))\right) + O_g\left(\frac{1}{l}\right) \frac{\#G(l^2)}{\#B(l^2)} x^{1/2} + x^{1/2}$$

$$\ll l^{2g^2 - 1} \log(l) + l^{2g^2 - 2g + 1} x^{1/2} + x^{1/2}.$$

We therefore have (since $l \leq x$ and $\theta \geq 1/2$)

$$\#\mathcal{A}_{l^2} \ll_g \frac{1}{l^2} \operatorname{li}(x) + l^{(9/2)g^2 - 1/2} x^\theta \log x + x^{1/2},$$

so that, integrating by parts,

$$\sum_{(x^\xi)^V \leq l < (x^\xi)^U} \#\mathcal{A}_{l^2} \ll x^{-\xi U} \operatorname{li}(x) + x^{\theta + \xi U \left((9/2)g^2 + 1/2\right)} \log x + x^{1/2 + \xi U}$$

We therefore are ensured of (9) as long as $\boxed{\xi U < 1}$ and $\boxed{\theta + \xi U \left((9/2)g^2 + 1/2\right) < 1.}$

## 5.3 Determining the optimal constants.

Collecting the constraints, we see that our goal is achieved as long as

$$g < \xi(rU + V), \quad \theta + \xi\left(g^2 + \frac{3}{2}g + 3\right) < 1, \quad \xi U < 1, \quad \theta + \xi U \left((9/2)g^2 + 1/2\right) < 1.$$

To attain minimal $r$, we minimize the value of
$$\frac{1}{U}\left(\frac{g}{\xi} - V\right),$$
so we wish to maximize $\xi$, $U$, and $V$ within our constraints.

Certainly, the constraint $\xi U < 1$ is redundant. Recall that the constraints of the sieve include $V \leq 1/4$ and $1/2 \leq U < 1$. We thus choose $V = 1/4$. Then, in particular, the terms $\alpha(V) = 0$ and $\beta(V) = 0$. Thus, doing a bit of calculus, we see that in order for $J(U, 1/4) > 0$, so that $B > 0$, we must have $\boxed{U < 3/4.}$

Thus, take
$$\xi = \frac{1-\theta}{(9/2)g^2 + 1/2}\left(\frac{4}{3} + \epsilon\right); \quad U = \frac{3}{4} - \epsilon.$$

Then, we see that for $g \geq 2$, the constraint $\theta + \xi\left(g^2 + \frac{3}{2}g + 3\right) < 1$ is satisfied for any $\epsilon > 0$. Thus, we may take
$$r = \left\lceil \frac{(9/2)g^3 + (1/2)g}{1 - \theta} - \frac{1}{3} \right\rceil$$
and $\epsilon$ sufficiently small. This concludes the proof of Theorem 1.3.

### 5.4 Proof of Theorem 1.5

We follow the argument of David-Wu to prove Theorem 1.5. Write the usual sieving function,
$$S(\mathcal{A}, \mathcal{P}, z) := \#\left(\mathcal{A} \setminus \bigcup_{p \in \mathcal{P}, p \leq z} \mathcal{A}_p\right).$$

Then, from the Weil bound, we see that for any $z < x$,
$$\#\left\{p \leq x \mid \#A_p(\mathbb{F}_p) \text{ is prime}\right\} = \#\left\{p \leq x \mid \#A_p(\mathbb{F}_p) \text{ is prime}, \#A_p(\mathbb{F}_p) > z\right\}$$
$$+ \#\left\{p \leq x \mid \#A_p(\mathbb{F}_p) \text{ is prime}, \#A_p(\mathbb{F}_p) \leq z\right\}$$
$$\leq S(\mathcal{A}, \mathcal{P}, z) + O(z^{1/g}).$$

We now apply the Selberg linear sieve (see Theorem 8.3 of [HR74]), with $q = 1$, and in their notation, $\xi = z$, which yields
$$S(\mathcal{A}, \mathcal{P}, z) \leq XV(z)\left(F(2) + o(1)\right) + \mathcal{R}$$
where
$$\mathcal{R} = \sum_{d < z^2, d \mid P(z)} 3^{\omega(d)} |r(\mathcal{A}, d)|$$
$$\ll_g \sum_{d < z^2} \mu(d)^2 3^{\omega(d)} d^{g^2 + (3/2)g + 1} x^\theta \log(x)$$
$$\ll x^\theta z^{2g^2 + 3g + 6} \log(x)$$
which is $o\left(x/(\log x)^2\right)$ if $\log(z)/\log(x) < (1-\theta)/(2g^2 + 3g + 6)$.

Choose $\epsilon > 0$ and define $z$ via $\log(x)/\log(z) = (2g^2 + 3g + 6)/(1 - \theta) + \epsilon$. Then, the definition of $F(u)$ tells us that $F(2) = e^\gamma$, so
$$X \cdot V(z)\left(F(2) + o(1)\right) = C_A \frac{\text{li}(x)}{\log z} \cdot \prod_{l > z}\left(\frac{1 - \#C(l)/\#G(l)}{1 - 1/l}\right)^{-1}(1 + o(1))$$
$$= C_A \frac{x}{(\log x)^2} \cdot \frac{\log(x)}{\log(z)} \cdot (1 + o(1))$$
$$\leq \left(\frac{2g^2 + 3g + 6}{1 - \theta} + \epsilon'\right) C_A \frac{x}{(\log x)^2}$$

for $x \gg_{A,\theta,\epsilon'} 0$, and the result follows.

## 5.5 Proof of Theorem 1.6

We continue the assumption that $A/\mathbb{Q}$ is generic and that the $\theta$-Hypothesis holds for $A$. We will employ Theorem 3.8 with the data

$$S := \{p \leq x\};$$
$$f(p) := \#A_p(\mathbb{F}_p);$$
$$\lambda_l := \#C(l)/\#G(l);$$

and the functions $e_l(x)$ and $e_{l_1\ldots l_u}(x)$ defined accordingly. We let $\beta \in (0,1]$ be arbitrary, and $\alpha = \alpha(x)$ arbitrary such that $0 < \alpha(x) < \beta$. We define $y = x^\alpha$, and will determine sufficient conditions on $\alpha$ and $\beta$ for conditions (1)-(6) of Theorem 3.8 to be satisfied.

We note that our choice of $S$ does not agree with our methods in this article so far; here, we do not exclude those $p$ for which $\#A_p(\mathbb{F}_p)$ shares a factor with $M_A$. It is clear, though, that the bound $r(\mathcal{A}, d) \ll d^{g^2+(3/2)g+1} x^\theta \log x$, for squarefree $d$, holds as well for the error function in this context: that is,

$$\pi(x) \cdot e_d(x) \ll d^{g^2+(3/2)g+1} x^\theta \log x.$$

We proceed:

1. Let $p \in S(x)$. Then, $f(p) = (1+o(1))p^g$, by the Weil Conjectures. Thus, for any chosen $\beta$, the number of distinct prime divisors of $f(p)$ that are more than $x^\beta$ is bounded by $(\log(g) + o(1))/\log(\beta)$.

2. We have

$$\sum_{y<l<x^\beta} \lambda_l = \sum_{y<l<x^\beta} l^{-1} + O_g(l^{-2})$$
$$= \log\log(x^\beta) - \log\log(x^\alpha) + O(1)$$
$$= -\log\alpha + O(1).$$

We must thus require $\log\alpha = o(\sqrt{\log\log x})$.

3. We have

$$\sum_{y<l<x^\beta} |e_l| \ll \sum_{y<l<x^\beta} l^{g^2+(3/2)g+1} x^{-1+\theta} (\log x)^2$$
$$\leq x^{-1+\theta} (\log x)^2 x^{\beta \cdot (g^2+(3/2)g+2)}$$

This quantity is $o(\sqrt{\log\log x})$ if

$$\beta < \frac{1-\theta}{g^2+(3/2)g+2},$$

and $\alpha$ satisfies the condition in item 2. (Since $\theta < 1$, such a $\beta$ exists.)

4. As in item 2, we have

$$\sum_{l \leq y} \lambda_l = \log\log x + \log\alpha + O(1)$$

which is of the desired form (with the constant $c = 1$) assuming the condition in item 2.

5. The quantity $\sum_{l \leq y} \lambda_l^2$ is clearly $O(1)$ from the previous discussion.

6. Lastly, mimicking [Liu06], we have

$$\sum_\star |e_{l_1 \cdot l_u}| \ll x^{-1+\theta} (\log x)^2 \left( \sum_{l \leq x^\alpha} l^{g^2+(3/2)g+1} \right)^u$$
$$\ll x^{-1+\theta+\alpha \cdot u \cdot (g^2+(3/2)g+2)} (\log x)^2$$

which, assuming that $\alpha(x) \to 0$, is asymptotic to $x^{-1+\theta+o(1)} = o((\log\log x)^{-r/2})$ for any $r$, since $\theta < 1$.

We thus require the existence of $\alpha(x)$ such that $\alpha = o(1)$ but $\log(\alpha(x)) = o(\sqrt{\log\log x})$, which is clear: take, for instance, $\alpha = (\log\log x)^{-1}$. This concludes the proof of Theorem 1.6.

# 6 A Koblitz Conjecture for Higher Genus and Experimental Evidence

The heuristics of the Koblitz Conjecture suggest the following conjecture.

**Conjecture 6.1.** *Let $A/\mathbb{Q}$ be an abelian variety satisfying the hypothesis $(\text{Triv}_A)$ such that $C_A \neq 0$. Then,*

1. $\#\{p \leq x \mid \#A_p(\mathbb{F}_p) \text{ is prime}\} \asymp_A \dfrac{x}{(\log(x))^2}.$

2. *In particular, if $A$ is generic,* $\#\{p \leq x \mid \#A_p(\mathbb{F}_p) \text{ is prime}\} \sim C_A \dfrac{x}{g(\log(x))^2}.$

Our Conjecture appears to be consistent with the generalizations by Weng and Spreckels, but we do not know of any author who has posited this Conjecture as stated (that is, for a *fixed* $A/\mathbb{Q}$). We also believe that part (2) of Conjecture 6.1 could be extended to those abelian varieties $A$ with $\text{End}(A)$ larger than $\mathbb{Z}$, analogously to Conjecture B of [Kob88], but we hesitate to do so for concern about stating the asymptotic constant correctly.

We provide experimental evidence for Conjecture 6.1 in the reminder of this Section. We collected from the LMFDB [LMFDB] some hyperelliptic curves $C/\mathbb{Q}$ of genus $g = 2$ whose Jacobians $J_C$ are generic and satisfy condition $(\text{Triv}_{J_C})$. We also considered the hyperelliptic genus 3 curve $C_3$ given by the equation

$$y^2 = x^7 - 14085x^6 + 33804x^5 - 27231x^4 + 27231x^3 - 35995x^2 - 33803x + 25039;$$

this curve was produced in the recent paper of Arias-de-Reyna et al. [Ari+16] as an example of a genus 3 curve whose Jacobian is proven to be generic by their Theorem 4.1. We ran a Sage program to collect the group orders $\#(J_C)_p(\mathbb{F}_p)$, with $p \leq 2^{20}$ for the genus 2 curves, and $p \leq 6 \cdot 10^4$ for $C_3$. (We had difficulty computing the group orders for larger $p$.) We then graphed the ratio

$$\frac{\#\left\{p \leq x \,\middle|\, \#(J_C)_p(\mathbb{F}_p) \text{ is prime}\right\}}{\pi(x)/(\log x)}$$

for $x$ at prime values $q$ for which $\#(J_C)_q(\mathbb{F}_q)$ is prime. We display these graphs in Figure 1. This evidence supports part (1) of the Conjecture, and if we were able to compute the constant $C_A$, we could check whether the evidence also supports part (2).

In the spirit of the questions of Lang-Trotter results "on average" (see, for instance, [BCD11]), we also approximated what we might call "universal constants"

$$\mathfrak{C}_g := \prod_\ell \left( \frac{1 - \#C(\ell)/\#G(\ell)}{1 - 1/\ell} \right)$$

where $G(l) = \text{GSp}_{2g}(\mathbb{Z}/l\mathbb{Z})$, and $C(l)$ is the union of conjugacy classes in $G(l)$ defined in (10). For a given generic abelian variety $A$, the constant $C_A$ differs from $\mathfrak{C}_g$ only by a factor depending on its non-surjective primes. We computed these approximations by finding the product for $\ell < 2^n$, for $n \leq 24$; they appear in Figure 2.

Interestingly, the functions for genus 2 curves in Figure 1 appear to converge to values which differ from $\mathfrak{C}_2/2$ by approximately half. This is perhaps more than one might expect: the author expects that the Euler factors at which $C_A$ and $\mathfrak{C}_2$ disagree (namely, those for the non-surjective primes of $A$) are not significantly different in magnitude, and he expects that there are not many such Euler factors.

It also appears that the limit $\lim_{g \to \infty} \mathfrak{C}_g$ exists. Very similar constants were computed in [Cas+12] in the context of Jacobians of hyperelliptic curves, though of course once $g \geq 3$, not all curves are hyperelliptic, and once $g \geq 4$, not all ppav's are Jacobians.

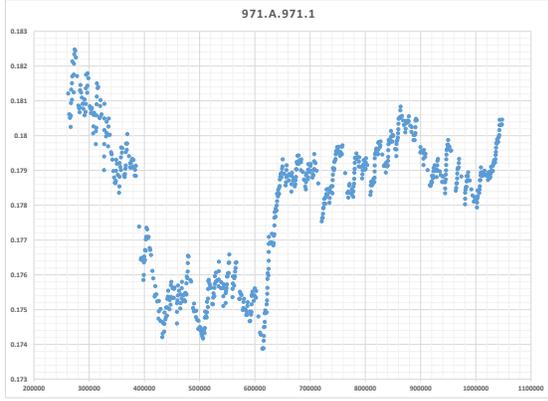

(a) curve 971.a.971.1 with equation
$Y^2 + Y = X^5 - 2X^3 + X$

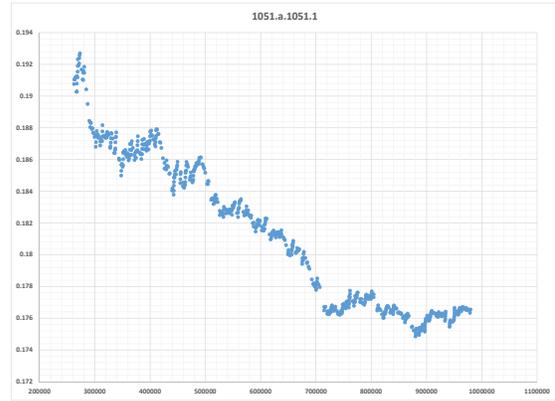

(b) curve 1051.a.1051.a with equation
$Y^2 + Y = X^5 - X^4 + X^2 - X$

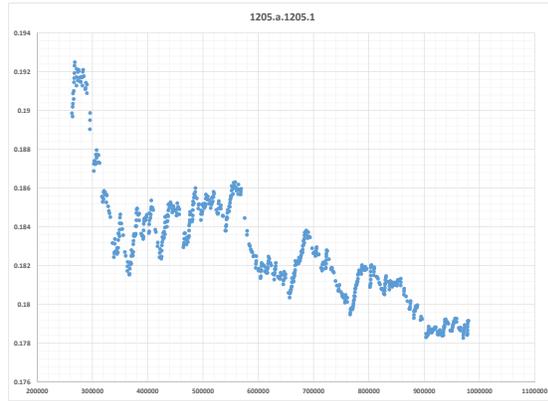

(c) curve 1205.a.1205.1 with equation
$Y^2 + Y = X^5 + 2X^4 - X^2$

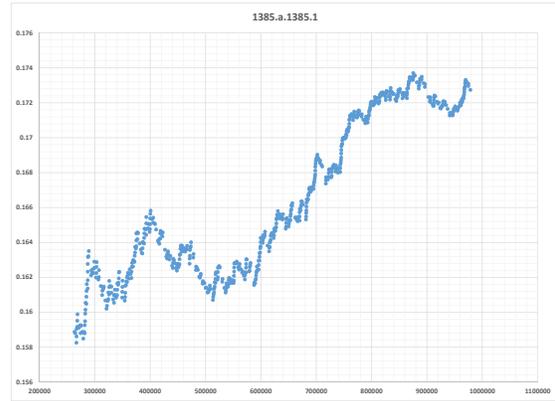

(d) curve 1385.a.1385.1, with equation
$Y^2 + Y = X^5 + 3X^4 + 3X^3 - X$

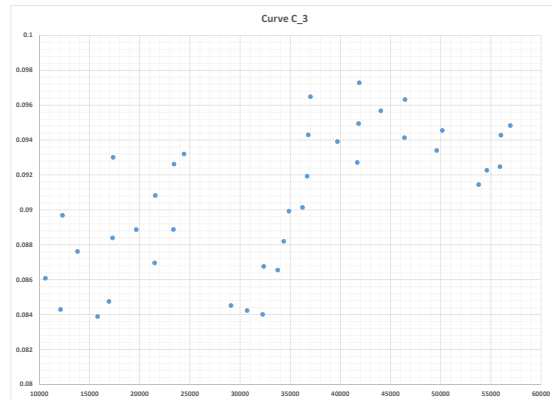

(e) genus 3 curve $C_3$

Figure 1: $\dfrac{\#\{p \leq x \mid \#(J_C)_p(\mathbb{F}_p) \text{ is prime}\}}{\pi(x)/(\log x)}$, with $x \in [2^{18}, 2^{20}]$ for genus 2 curves, $x \in [10^4, 6 \cdot 10^4]$ for $C_3$

| $n$ | $\mathfrak{C}_1$ | $\mathfrak{C}_2$ | $\mathfrak{C}_3$ | $\mathfrak{C}_4$ |
|-----|------------------|------------------|------------------|------------------|
| 2   | 0.562500000000000 | 0.760989583333333 | 0.754354887320847 | 0.754413616554689 |
| 4   | 0.513926644244210 | 0.706235456622878 | 0.700012977803311 | 0.700067571267533 |
| 8   | 0.505468861944026 | 0.695053638628807 | 0.688929626754209 | 0.688983355837062 |
| 16  | 0.505166809270517 | 0.694639169901420 | 0.688521872595408 | 0.688572506891267 |
| 24  | 0.505166169952616 | 0.694638290801478 | 0.688517938493554 | 0.688571635469346 |

Figure 2: Computations for the constants $\mathfrak{C}_g$.